\documentclass[12pt]{article}
\usepackage{mathrsfs}
\usepackage{graphicx,amsmath,amsfonts,amssymb,amscd,amsthm}
\newtheoremstyle{kai}
{3pt} {3pt} {} {} {\bfseries} {.} {.5em} {}
\def\EquationsBySection{\def\theequation
{\thesection.\arabic{equation}}%
\@addtoreset{equation}{section}}

\newcommand\old[1]{}
\newcommand{\pend}{\hfill \thicklines \framebox(6.6,6.6)[l]{}}

\renewenvironment{proof}{\noindent {\it  Proof.} \rm}{\pend}
\newtheorem{theorem}{Theorem}[section]
\newtheorem{lemma}{Lemma}[section]
\newtheorem{corollary}{Corollary}[section]

\newtheorem{remark}{Remark}[section]

\newtheorem{example}{Example}[section]

\EquationsBySection \makeatother
\usepackage{indentfirst}

\newcommand\RR{\mathbb{R}}

\begin{document}
\pagestyle{plain}
\title
{\bf Stochastic Population Dynamics Driven by L\'{e}vy Noise}
\author{Jianhai Bao and  Chenggui Yuan\thanks{{\it E-mail address:}
C.Yuan@swansea.ac.uk.}\\
Department of Mathematics,\\
 Swansea University, Swansea SA2 8PP, UK
\\
}
\date{}
\maketitle
\begin{abstract}{\rm This paper considers stochastic population dynamics driven by
L\'{e}vy noise. The contributions of this paper lie in that (a)
Using Khasminskii-Mao theorem, we show that the stochastic
differential equation associated with the model has a unique global
positive solution; (b) Applying an exponential martingale inequality
with jumps,  we discuss the
asymptotic pathwise estimation of such model. }\\

\noindent {\bf Keywords:} Brownian motion, L\'{e}vy noise; Exponential martingale
inequality with jumps.\\
\noindent{\bf Mathematics Subject Classification (2000)} \ 93D05,
60J60, 60J05.
\end{abstract}

 \noindent

\section{Introduction}
Stochastic population dynamics perturbed by Brownian motion has been
studied extensively by many authors. There are a great amount of
literature on this topic. In \cite{mmr02} Mao, Marion and Renshaw
investigate
 stochastic $n$-dimensional Lotka-Volterra system
\begin{equation}\label{eq126}
dX(t)=\mbox{diag}(X_1(t),\cdots,X_n(t))\left[(b+AX(t))dt+\sigma
X(t)dW(t)\right],
\end{equation}
where
\begin{equation*}
X=(X_1,\cdots,X_n)^T,\ \ b=(b_1,\cdots,b_n)^T, \ \
A=(a_{ij})_{n\times n},\  \ \sigma=(\sigma_{ij})_{n\times n},
\end{equation*}
and reveal the important role that the environmental noise can
suppress a potential population explosion; Mao, Sabanis and Renshaw
\cite{msr03} further discuss the asymptotic behaviour of
 population process determined by Eq.
\eqref{eq126}. Then the techniques developed in \cite{mmr02,msr03}
have been applied successfully to study stochastic delay population
dynamics in \cite{bm04,myz05}; functional Kolmogorov-type systems
\cite{wh08,wh10}; hybrid competitive Lotka-Volterra models in
\cite{my06,yml09,zy09a,zy09b}.

The population may suffer sudden environmental shocks, e.g.,
earthquakes, hurricanes, epidemics, etc. However, stochastic
extension of population process described by Eq. \eqref{eq126}
cannot explain the phenomena above. To explain these phenomena,
introducing a jump process into underlying population dynamics is
one of the important methods. To our knowledge there is few
systematic work so far in which the noise source is a jump process.
The work presented here is to take some steps in this direction,
building extensively on  the existing results mentioned above for
the Brownian motion case,  is also a sequel to that of
\cite{bmyy11}, where competitive Lotka-Volterra population dynamics
with jumps
\begin{equation*}
dY(t)=Y(t)\Big[\Big(a(t)-b(t)Y(t)\Big)dt+\sigma(t)dW(t)
+\int_{\mathbb{Y}}\gamma(t,u)\tilde{N}(dt,du)\Big]
\end{equation*}
is investigated, and  the explicit solution, sample Lyapunov
exponent and invariant measure are also addressed.

We focus in this paper on stochastic population dynamics
\eqref{eq126} that is further perturbed by L\'{e}vy noise, that is,
\begin{equation}\label{eq1}
\begin{split}
dX(t)&=\mbox{diag}(X_1(t),\cdots,X_n(t))\Big[(b+AX(t))dt+\sigma
X(t)dW(t)\\
&\ \ \ \ + \int_{\mathbb{Y}}H(X(t^-),u)\tilde{N}(dt,du)\Big].
\end{split}
\end{equation}
Here $b,A,\sigma$ are defined as in model \eqref{eq126}, $W(t)$ is a
scalar Brownian motion defined on the probability space
$\{\Omega,{\mathcal F}, \mathbb{P}\}$ with the filtration
$\{{\mathcal F}_{t}\}_{t\geq0}$ satisfying the usual condition,
$X(t^-):=\lim_{s\uparrow t}X(s)$,  $N(dt,du)$ is a real-valued
Poisson counting measure with characteristic measure $\lambda$ on
measurable subset $\mathbb{Y}$ of $[0,\infty)$ with
$\lambda(\mathbb{Y})<\infty$,
$\tilde{N}(dt,du):=N(dt,du)-\lambda(du)dt$, and
$H:\mathbb{R}^n\times\mathbb{Y}\rightarrow\mathbb{R}^n$. Throughout
this paper, we further assume that $W$ is independent of $N$.

In reference to the existing results in the literature, our
contributions are as follows:
\begin{itemize}
\item We use jump diffusion to model the evolutions of
population dynamics when they suffer sudden environmental shocks;
\item Using  Khasminskii-Mao theorem,  we
show that the Stochastic Differential Equation (SDE) associated with
the model has a unique global positive solution; \item Applying an
exponential martingale inequality with jumps,  together with the
standard Borel-Cantelli lemma, we discuss the asymptotic pathwise
estimation of such model.
\end{itemize}

\section{Global Positive Solutions}

Since $X$ denotes the population sizes of the $n$ interacting
species, it is natural to require the solution of Eq. \eqref{eq1}
not only to be positive but also not to explode in a finite time.
Therefore, in this section we intend to show that Eq. \eqref{eq1}
has a unique global positive solution under some conditions. Since
the coefficients don't satisfy linear growth condition or weak
coercivity condition,  even they satisfy local Lipschitz condition,
the solutions of Eq. \eqref{eq1} may explode in a finite time.
Khasminskii \cite[Theorem 4.1, p85]{k80} and Mao \cite{Mao02} gave
the Lyapunov function argument, which is a powerful test for
nonexplosion of solutions without linear growth contion and is
referred as Khasminskii-Mao theorem. In what follows, we shall also
apply  Khasminskii-Mao approaches to show that Eq. \eqref{eq1} has a
unique global positive solution $X(t), t\geq0$. In this section we
will show:

\begin{itemize}
\item  Jump processes can suppresses the explosion;

\item   Brownian motion can also suppresses the explosion to our new model, which is similar to that of  \cite{mmr02}.
\end{itemize}
 In what follows let $K>0$ be a
generic constant whose values may vary for its different
appearances. To show the main result let us recall the following
facts.

Consider $1$-dimensional SDE with jumps
\begin{equation}\label{eq00}
dX(t)=F(X(t))dt+G(X(t))dW(t)+\int_{\mathbb{Y}}\Phi(X(t^-),u)\tilde{N}(dt,du),\
\  t\geq0
\end{equation}
with initial condition $X(0)=x_0\in\mathbb{R}$, where $W(t)$ is a
real-valued  Brownian motion, $F,G:\mathbb{R}\rightarrow\mathbb{R},
\Phi:\mathbb{R}\times\mathbb{Y}\rightarrow\mathbb{R}$.

The following conclusion is given by  (see \cite[Lemma 2]{wee}.
\begin{lemma}\label{lemma}
{\rm Let $F(0)=G(0)=\Phi(0,u)=0$ for $u\in\mathbb{Y}$ and $F,G,\Phi$
satisfy local Lipschitz condition. Set
\begin{equation*}
J(x):=\int_{\mathbb{Y}}\left(\ln\frac{|x+\Phi(x,u)|}{|x|}\right)^2\lambda(du), \, x \neq 0.
\end{equation*}
Assume  that
\begin{equation}\label{eq103}
\sup_{0<|x|\leq m}J(x)<\infty \ \ \ \mbox{ for each } m\ge 1.
\end{equation}
Then for $x\neq0$
\begin{equation*}
\mathbb{P}(X(t,x)\neq0 \mbox{ and } X(t^-,x)\neq0\mbox{ for any }
t\geq0)=1,
\end{equation*}
where $X(t,x)$ denotes the solution of Eq. \eqref{eq00} starting
from $x$ at time $t=0$. }
\end{lemma}

For convenience of reference, we recall some fundamental
inequalities stated as a lemma.
\begin{lemma}
{\rm \begin{equation}
\label{eq108} x^r\leq1+r(x-1),\ \ \ \ x\geq0,\
\ \ 1\geq r\geq0,
\end{equation}
\begin{equation}\label{eq109}
n^{(1-\frac{p}{2})\wedge0}|x|^{p}\leq\sum_{i=1}^nx_i^p\leq
n^{(1-\frac{p}{2})\vee0}|x|^{p}, \ \ \forall p>0,\ \
x\in\mathbb{R}^n_+,
\end{equation}
where $\mathbb{R}^n_+:=\{x\in \mathbb{R}^n:x_i>0, 1\leq i\leq n\}$,
and
\begin{equation}\label{eq110}
\ln x\leq x-1,\ \ \ \ \  \ \ \ \ \ \ \ \ \ \ \ \  x>0.
\end{equation}

}

\end{lemma}

\subsection{Explosive Suppression by Jump Processes}

For jump-diffusion coefficient we  assume that

\noindent(${\bf H1}$) For any $m\geq1,
x\in\mathbb{R}^n,u\in\mathbb{Y}$ and $i=1,\cdots, n$
\begin{equation}\label{eq114}
H_i(x,u)>-1, \ \ \ \ \ H_i(0,u)=0,
\end{equation}
\begin{equation}\label{eq105}
\sup_{0<|x|\leq
m}\int_{\mathbb{Y}}(\ln|1+H_i(x,u)|)^2\lambda(du)<\infty,
\end{equation}
and  for each $k>0$ there exists constant $L_k>0$ such that
\begin{equation}\label{eq104}
\int_{\mathbb{Y}}|H(x,u)-H(y,u)|^2\lambda(du)\leq L_k|x-y|^2
\end{equation}
whenever $x,y\in \mathbb{R}^n$ with $|x|\vee|y|\leq k$.

\begin{theorem}\label{general2}
{\rm Let assumption $({\bf H1})$ hold. Assume further that for
$p\in(0,1)$ there exist constants $\delta>0,\alpha>2$ such that for
$x\in\mathbb{R}^n, i=1,\ldots n,$
\begin{equation}\label{eq11}
J_i(x,p):=\int_{\mathbb{Y}}\left[(1+H_i(x,u))^{p}-1
-pH_i(x,u)\right]\lambda(du)\leq-\delta|x|^\alpha+\mbox{o}(|x|^\alpha),
\end{equation}
where $\mbox{o}(|x|^\alpha)/|x|^\alpha\rightarrow0$ as
$|x|\rightarrow\infty$. Then, for any initial condition
$\bar{x}\in\mathbb{R}^n_+$, Eq. \eqref{eq1} has a unique global
solution $X(t)\in\mathbb{R}^n_+$ for any $t\geq0$ almost surely.}
\end{theorem}

Before we prove the theorem, we give an example such that condition \eqref{eq11} holds.

\begin{example}\label{example}
{\rm For $i=1,\cdots,n$, let $\gamma_i>0$ and
$\int_\mathbb{Y}(1\vee\gamma_i(u))\lambda(du)<\infty$. Assume that
\begin{equation*}
H_i(x,u):=\gamma_i(u)H(|x|), \ \ \ x\in\mathbb{R}^n, u\in\mathbb{Y},
\end{equation*}
where $H(|x|)$ is a polynomial of degree $\alpha>2$ with positive
leading coefficient. Then by a straightforward computation we have
for some $\delta>0$
\begin{equation*}
\begin{split}
J_i(x,p)&=\int_{\mathbb{Y}}\left[(1+\gamma_i(u)H(|x|))^{p}-1
-p\gamma_i(u)H(|x|)\right]\lambda(du)\\
&\leq-\delta|x|^\alpha+\mbox{o}(|x|^\alpha).
\end{split}
\end{equation*}
Therefore condition \eqref{eq11} holds. }
\end{example}

\noindent{\bf Proof of Theorem \ref{general2}.}
 By \eqref{eq104}, for arbitrary initial value $\bar{x}\in\mathbb{R}^n_+$ there is a unique local solution
 $X(t)$ for $t\in[0,\tau_e)$, where $\tau_e$ is the explosion time. By Eq. \eqref{eq1} the $i$th component
$X_i(t)$ of $X(t)$  admits the form
\begin{equation}\label{eq12}
\begin{split}
dX_i(t)&=X_i(t)\Big[\Big(b_i+\sum\limits_{j=1}^na_{ij}X_j(t)\Big)dt+\sum\limits_{j=1}^n\sigma_{ij}X_j(t)dW(t)\\
&\quad\ \ \ \ \ \ \ \ \
+\int_{\mathbb{Y}}H_i(X(t^-),u)\tilde{N}(dt,du)\Big].
\end{split}
\end{equation}
 Note that for any $t\in[0,\tau_e)$
\begin{equation*}
\begin{split}
X_i(t)=\bar{x}_i\exp\Big\{&\int_0^t\Big(b_i+\sum\limits_{j=1}^na_{ij}X_j(s)-\frac{1}{2}
\Big(\sum\limits_{j=1}^n\sigma_{ij}X_j(s)\Big)^2\\
&+\int_\mathbb{Y}(\ln(1+H_i(X(s),u))-H_i(X(s),u))\lambda(du)\Big)ds\\
&+\int_0^t\sum\limits_{j=1}^n\sigma_{ij}X_j(s)dW(s)+\int_0^t\int_\mathbb{Y}\ln(1+H_i(X(s^-),u))\tilde{N}(ds,du)\Big\}.
\end{split}
\end{equation*}
This, together with $\bar{x}\in\mathbb{R}^n_+$, yields that
$X_i(t)\geq0$ for any $t\in[0,\tau_e)$. On the other hand, due to
\eqref{eq105}, for Eq. \eqref{eq12} condition \eqref{eq103} holds.
Then Lemma \ref{lemma} gives $X_i(t)>0$ for any $t\in[0,\tau_e)$
since $\bar{x}\in\mathbb{R}^n_+$. Next we show $\tau_e=\infty$ a.s.
Let $k_0>0$ be sufficiently large such that $|\bar{x}|< k_0$. For
each $k> k_0$ define a stopping time
\begin{equation*}
\tau_k:=\inf\{t\in[0,\tau_e):|X(t)|>k \}.
\end{equation*}
Clearly, $\tau_k$ is increasing as $k\uparrow\infty$. Set
$\tau_\infty:=\lim_{k\rightarrow\infty}\tau_k$, whence
$\tau_\infty\leq\tau_e$ a.s., then it is sufficient to check
$\tau_{\infty}=\infty$ a.s. Introduce a Lyapunov function for any
$p\in(0,1)$
\begin{equation}\label{eq107}
V(x):=\sum_{i=1}^nx_i^p,\ \  x\in\mathbb{R}^n_+.
\end{equation}
Let $T>0$ be arbitrary. For any $0\leq t\leq  \tau_k\wedge T$, the
It\^o formula yields
\begin{equation}\label{eq111}
\begin{split}
dV(X(t))&=\mathcal {L}V(X(t))dt+p\sum\limits_{i=1}^nX_i^p(t)\sum\limits_{j=1}^n\sigma_{ij}X_j(t)dW(t)\\
&\quad+\sum\limits_{i=1}^n\int_{\mathbb{Y}}[(1+H_i(X(t^-),u))^p-1]\tilde{N}(dt,du)X_i^p(t),
\end{split}
\end{equation}
where for $x\in\mathbb{R}^n_+$
\begin{equation}\label{eq112}
\begin{split}
\mathcal {L}V(x)&:=p\sum\limits_{i=1}^n
\left[b_i+\sum\limits_{j=1}^na_{ij}x_j-\frac{1-p}{2}\left(\sum\limits_{j=1}^n\sigma_{ij}x_j\right)^2\right]x_i^{p}
\\
&\quad+\sum\limits_{i=1}^n\int_{\mathbb{Y}}\left[(1+H_i(x,u))^{p}-1
-pH_i(x,u)\right]\lambda(du)x_i^{p}\\
&:=K_1(x,p)+K_2(x,p).
\end{split}
\end{equation}
By inequality \eqref{eq109}
\begin{equation*}
K_1(x,p)\leq K|x|^{2+p}+\mbox{o}(|x|^{2+p}),
\end{equation*}
and, thanks to \eqref{eq11}
\begin{equation*}
K_2(x,p)\leq-\delta|x|^{\alpha+p}+\mbox{o}(|x|^{\alpha+p}).
\end{equation*}
Thus, for $\alpha>2$
\begin{equation}\label{eq78}
\mathcal {L}V(x)\leq K \ \ \mbox{ for any } x\in\mathbb{R}^n_+.
\end{equation}
Define for each $u>0$
\begin{equation*}
\mu(u):=\inf\{V(x), |x|\geq u\}.
\end{equation*}
Thanks to inequality \eqref{eq109}, it is easy to see that
\begin{equation}\label{eq22}
\lim\limits_{u\rightarrow\infty}\mu(u)=\infty.
\end{equation}
 Then we obtain from \eqref{eq78} that for some constant $K>0$
\begin{equation*}
\mu(k)\mathbb{P}(\tau_k\leq
T)\leq\mathbb{E}(V(X(\tau_k))I_{\{\tau_k\leq T\}})\leq
\mathbb{E}V(X(\tau_k\wedge T))\leq K.
\end{equation*}
Recalling \eqref{eq22} and letting $k\rightarrow\infty$ yields
\begin{equation*}
\mathbb{P}(\tau_{\infty}\leq T)=0.
\end{equation*}
Since $T$ is arbitrary, we must have
\begin{equation*}
\mathbb{P}(\tau_{\infty}=\infty)=1
\end{equation*}
and Eq. \eqref{eq1} admits a unique global solution
$X(t)\in\mathbb{R}^n_+$ on $t\geq0$.

\begin{remark}
{\rm In \cite{mmr02}, under the condition

\noindent(${\bf H2}$) \ \ \ \ \ \ \ \ \ \ \ \ \ \ \ \
$\sigma_{ii}>0$ if $1\leq i\leq n$ while $\sigma_{ij}\geq0$ if
$i\neq j$.\\ Mao, Marion and Renshaw reveal the important fact that
Brownian motion noise  can suppress a potential population
explosion. Theorem \ref{general2} shows that L\'{e}vy noise can also
play the same role, without any conditions being imposed on the
diffusion coefficient $\sigma$. }
\end{remark}

 As for population dynamics, in general, the following Lyapunov
function for $p\in(0,1)$
\begin{equation}\label{eq113}
U(x):=\sum\limits_{i=1}^n[x_i^p-1-p\ln x_i], \ \ \
x\in\mathbb{R}^n_+,
\end{equation}
is constructed to show that the SDE associated to the model admits a
unique global positive solution, see, e.g.,
\cite{mmr02,myz05,yml09,zy09a,zy09b}. In what follows, under
suitable conditions we can also show that Eq. \eqref{eq1} has a
unique global positive solution through the Lyapunov function
defined by \eqref{eq113}.

\begin{theorem}\label{existence}
{\rm Suppose that assumptions \eqref{eq114}, \eqref{eq104} and
\eqref{eq11} hold. Assume further that there exist constants
$\beta\in(0,\alpha]$ and $\nu>0$ such that
\begin{equation}\label{eq115}
\int_{\mathbb{Y}}\left[
H_i(x,u)-\ln(1+H_i(x,u))\right]\lambda(du)\leq
\nu|x|^\beta+\mbox{o}(|x|^{\beta})
\end{equation}
for $i=1,\cdots, n$ and $x\in\mathbb{R}^n_+$. Then, for any initial
condition $\bar{x}\in\mathbb{R}^n_+$, Eq. \eqref{eq1} has a unique
global solution $X(t)\in\mathbb{R}^n_+$ for any $t\geq0$ almost
surely. }
\end{theorem}

\begin{proof}
 Since the argument is similar to that of \cite[Theorem
2.1]{mmr02}, we here only sketch the proof to point out the
variation from the Brownian motion case. Let $k_0\in\mathbb{N}$ be
sufficiently large
 such that every component of $\bar{x}$ is contained in the interval
 $(\frac{1}{k_0},k_0)$. For each $k>k_0$ define a stopping time
\begin{equation*}
\tau_k:=\inf\Big\{t\in[0,\tau_e):X_i(t)\notin\Big(\frac{1}{k},k\Big)
\mbox{ for some } i=1,\cdots,n\Big\},
\end{equation*}
where $\tau_e$ is the explosion time. In the sequel, we show
$\tau_{\infty}:=\lim_{k\rightarrow\infty}\tau_k=\infty$ a.s. Let
$T>0$ be arbitrary. For any $0\leq t\leq  \tau_k\wedge T$, applying
It\^o's formula, we obtain
\begin{equation}\label{eq23}
\begin{split}
\mathcal {L}U(x):&=p\sum\limits_{i=1}^n\left[b_i(x_i^p-1)
-(x_i^p-1)\sum\limits_{j=1}^na_{ij}x_j
+\left(\frac{p-1}{2}x_i^p+1\right)\left(\sum\limits_{j=1}^n\sigma_{ij}x_j\right)^2\right]\\
&\quad+\sum\limits_{i=1}^n\int_{\mathbb{Y}}[(1+H_i(x,u))^p-1-pH_i(x,u)]\lambda(du)x_i^p\\
&\quad+p\sum\limits_{i=1}^n\int_{\mathbb{Y}}\left[
H_i(x,u)-\ln(1+H_i(x,u))\right]\lambda(du)\\
&:=I_1(x,p)+I_2(x,p)+I_3(x,p).
\end{split}
\end{equation}
By \eqref{eq109}, note that
\begin{equation*}
I_1(x,p)\leq K|x|^{2+p}+\mbox{o}(|x|^{2+p}).
\end{equation*}
Also, due to \eqref{eq11} and \eqref{eq115}
\begin{equation*}
I_2(x,p)+I_3(x,p)\leq
-\delta|x|^{\alpha+p}+np\nu|x|^\beta+\mbox{o}(|x|^{\alpha+p}).
\end{equation*}
 Then for $p\in(0,1)$ and $\beta\in(0,\alpha]$
\begin{equation*}
\mathcal {L}U(x)\leq K,\ \ x\in\mathbb{R}^n_+.
\end{equation*}
Define for each $u>1$
\begin{equation*}
\mu(u):=\inf\left\{U(x):  x_i\geq u \mbox{ or } x_i\leq\frac{1}{u}
\mbox{ for some } i=1,\cdots,n\right\}.
\end{equation*}
Due to the property of function $h(x):=x-1-\ln x, x>0$, we see that
\begin{equation*}
\lim\limits_{x\uparrow\infty}h(x)=\infty \mbox{ and }
\lim\limits_{x\downarrow0}h(x)=\infty
\end{equation*}
and hence
\begin{equation}\label{eq013}
\lim\limits_{u\rightarrow\infty}\mu(u)=\infty.
\end{equation}
The proof is then complete by repeating the procedure of
\cite[Theorem 2.1]{mmr02}.
\end{proof}

\subsection{Explosive Suppression  by Brownian Motion}
In this subsection, we further show that Brownian motion can also
suppresses the explosion to our model under condition $({\bf H2})$,
but weaker conditions imposed on jump-diffusion coefficient.

\begin{theorem}\label{general}
{\rm Assume that assumptions $({\bf H1})$ and $({\bf H2})$ hold.
Then, for any initial condition $\bar{x}\in\mathbb{R}^n_+$, Eq.
\eqref{eq1} has a unique global solution $X(t)\in\mathbb{R}^n_+$ for
any $t\geq0$ almost surely.}
\end{theorem}

\begin{proof}
Since the proof is very similar to that of Theorem \ref{general2},
we here only give an outline of the argument.  In \eqref{eq112},
note from $({\bf H2})$ and inequality \eqref{eq109} that for
$p\in(0,1)$ and $x\in\mathbb{R}^n_+$
\begin{equation}\label{eq24}
K_1(x,p)\leq -\frac{p(1-p)n^{-\frac{p}{2}}}{2}\min_{1\leq i\leq
n}\sigma_{ii} |x|^{2+p}+\mbox{o}(|x|^{2+p}).
\end{equation}
On the other hand, by inequality \eqref{eq108}, for any $p\in(0,1)$
and $x\in\mathbb{R}^n_+$
\begin{equation*}
\int_{\mathbb{Y}}\left[(1+H_i(x,u))^{p}-1
-pH_i(x,u)\right]\lambda(du)\leq0,
\end{equation*}
hence we have
\begin{equation*}
K_2(x,p)\leq0, \ \ \ p\in(0,1) \mbox{ and } x\in\mathbb{R}^n_+.
\end{equation*}
Thus
\begin{equation*}
\mathcal {L}V(x)\leq K \ \ \mbox{ for any } x\in\mathbb{R}^n_+.
\end{equation*}
The conclusion then follows by carrying out the procedure of Theorem
\ref{general2}.
\end{proof}

Applying Lyapunov function $U(x)$ in \eqref{eq113}, under suitable
conditions we can still guarantee that Eq. \eqref{eq1} admits a
unique global positive solution.

\begin{theorem}\label{theorem}
{\rm Let conditions \eqref{eq114}, \eqref{eq104}, $({\bf H2})$ hold
and assume further that condition \eqref{eq115} holds with
$\beta\in(0,2]$. Then, for any initial condition
$\bar{x}\in\mathbb{R}^n_+$, Eq. \eqref{eq1} has a unique global
solution $X(t)\in\mathbb{R}^n_+$ for any $t\geq0$ almost surely. }
\end{theorem}

\begin{proof}
The argument is similar to that of Theorem \ref{existence}. In Eq.
\eqref{eq23}, by inequality \eqref{eq108}
\begin{equation*}
I_2(x,p)\leq0, \ \ \ p\in(0,1) \mbox{ and } x\in\mathbb{R}^n_+.
\end{equation*}
What's more, condition \eqref{eq115} leads to,  for $\beta\in(0,2]$
and $x\in\mathbb{R}^n_+$,
\begin{equation*}
I_3(x,p)\leq np\nu|x|^\beta+\mbox{o}(|x|^{\beta}).
\end{equation*}
Combining \eqref{eq24},  for $p\in(0,1)$ and $\beta\in(0,2]$, we can
conclude that for some $K>0$
\begin{equation*}
\mathcal {L}U(x)\leq K,\ \ x\in\mathbb{R}^n_+.
\end{equation*}
Then the proof can be done by carrying out the procedure of Theorem
\ref{existence}.
\end{proof}

\begin{remark}
{\rm  By constructing two different Lyapunov functions $V(x)$ and
$U(x)$ defined by \eqref{eq107} and \eqref{eq113}, respectively,
under  different conditions, namely, \eqref{eq105} and
\eqref{eq115}, we show that Eq. \eqref{eq1} has a unique global
positive solution. Comparing conditions \eqref{eq105} and
\eqref{eq115}, we see that condition \eqref{eq105} is a local one
while   \eqref{eq115} has a growth restriction on jump diffusion $H$
in whole space. Therefore, the argument developed in Theorem
\ref{general} is easier to verify  than that of Theorem
\ref{theorem}.

}
\end{remark}

\section{Asymptotic Moment Properties}
In the last section, under suitable conditions we have shown that
Eq. \eqref{eq1} admits a unique global positive solution. From the
biological point of view, the nonexplosion property and positive
solution in a population dynamical system are often not good enough
while the moment properties  are more desired. In this section we
shall show that the $p$th moment with $p\in(0,1)$ and the average in
time of the moment of the solution to Eq. \eqref{eq1} are both
bounded.

\begin{theorem}\label{pth moment}
{\rm Under conditions of Theorem \ref{general2} (or Theorem
\ref{existence},  or Theorem  \ref{general}, or Theorem
\ref{theorem}), for any $p\in(0,1)$ and some $K>0$
\begin{equation}
\limsup_{t\rightarrow\infty}\mathbb{E}|X(t)|^p\leq K.
\end{equation}
}
\end{theorem}

\begin{proof}  We only give a sketch of proof under conditions of Theorem
\ref{general2} due to the similarities of arguments.
 For any
$p\in(0,1)$ let $V$ be defined by \eqref{eq107}.   For any
$|\bar{x}|<k$ define a stopping time
\begin{equation*}
\sigma_k:=\inf\{t\geq0:|X(t)|>k\}.
\end{equation*}
By the argument of Theorem \ref{general2}, we have $\sigma_k<\infty$
and $\sigma_k\uparrow\infty$ a.s. as $k\rightarrow\infty$. Applying
the It\^o formula yields
\begin{equation}
\mathbb{E}(e^{t\wedge\sigma_k}V(X(t\wedge\sigma_k)))=V(\xi(0))+\mathbb{E}\int_0^{t\wedge\sigma_k}e^{
s}[ V(X(s))+\mathcal {L}V(X(s))]ds,
\end{equation}
where $\mathcal {L}V$ is defined by \eqref{eq112}.
Since the leading terms of $\mathcal {L}V$ are $-\delta |x|^{\alpha+p}$ with $\alpha >2$
 we can deduce that there exists a constant $K>0$ such that $V(x)+\mathcal {L}V(x)\le K$. Hence
\begin{equation*}
\begin{split}
\mathbb{E}(e^{ t}V(X(t)))\leq K(1+e^{ t})
\end{split}
\end{equation*}
and the desired conclusion follows from inequality \eqref{eq109}.
\end{proof}

Going through the arguments of Theorem \ref{general2} (or Theorem
\ref{existence},  or Theorem  \ref{general}, or Theorem
\ref{theorem}), we can also derive that the average in time of the
moment of the solution to Eq. \eqref{eq1} is bounded.

\begin{theorem}
{\rm (i). Under conditions of Theorem \ref{general2} or Theorem
\ref{existence}, for any $p\in(0,1)$  and any initial value
$\bar{x}\in \mathbb{R}^n_+$, there exists constant $K>0$ such that
\begin{equation*}
\limsup_{t\rightarrow\infty}\frac{1}{t}\int_0^t\mathbb{E}|X(s)|^{p+2}ds\leq
K.
\end{equation*}

(ii). Under conditions of Theorem \ref{general} or Theorem
\ref{theorem}, for any $p\in(0,1)$ and any initial value $\bar{x}\in
\mathbb{R}^n_+$, there exists constant $K>0$ such that
\begin{equation*}
\limsup_{t\rightarrow\infty}\frac{1}{t}\int_0^t\mathbb{E}|X(s)|^{p+2}ds\leq
K.
\end{equation*}
}
\end{theorem}

\begin{proof}
Since the proofs of (i) and (ii) are very similar, we here only
prove (ii) under the conditions of Theorem \ref{general}. The
argument is motivated by that of  \cite[Theorem 2]{msr03}.
 By \eqref{eq111} and \eqref{eq112} we have
\begin{align*}
V(X(t))&\le V(X(0))+\int_0^t\bigg(p\sum_{i=1}^nb_i+p\sum_{i=1}^n\sum_{j=1}^nX_i^p(s)X_j(s)-\frac{p(p-1)}{2}\sum_{i=1}^n\sigma_{ii}^2X_i^{2+p}(s)\bigg)ds\\
&+M_1(t)+M_2(t),
\end{align*}
where $M_1(t), M_2(t)$ are two local martingales.  Noting that the polynomial
$$
p\sum_{i=1}^nb_i+p\sum_{i=1}^n\sum_{j=1}^nx_i^px_j-\frac{p(p-1)}{4}\sum_{i=1}^n\sigma_{ii}^2x_i^{2+p}
$$
has a upper bound $K$ (dependent on $p$), therefore
$$
V(X(t))+\frac{p(p-1)}{4}\sum_{i=1}^n\sigma_{ii}^2X_i^{2+p}(t)\le V(X(0))+Kt +M_1(t)+M_2(t).
$$
Taking the expectation,  and then dividing by $t$ on both sides, we obtain the result.
\end{proof}

By taking another different Lyapunov function, we will have the following theorem.

\begin{theorem}\label{theorem2}
{\rm Under the  conditions of Theorem \ref{general2},
 let $p^T=(p_1, \ldots p_n)$ be positive numbers such that $p_1+\ldots+p_n<1,$  and assume there exist constants $\beta_1$ and $\beta_2\in (0, \alpha)$ such that
 \begin{align}\label{eq11a}
 \int_{\mathbb{Y}}\bigg[ \Pi_{i=1}^n (1+H_i(x, u))^{p_i}-\sum_{i=1}^n(1+H_i(x, u))^{p_i}\bigg]\lambda(du)\le
 \beta_1|x|^{\beta_2}+o(|x|^{\beta_2}),
 \end{align}
where constant $\alpha>2$ was given in \eqref{eq11}. Then
$$
\mathbb{E} \left(\Pi_{i=1}^nX_i^{p_i}(t)\right)<\infty, \mbox{ for all } t\ge 0.
$$
 }
\end{theorem}

\begin{remark}
{\rm For $n=1$, condition \eqref{eq11a} must be true. Moreover,
Example \ref{example} also demonstrates that condition \eqref{eq11a}
holds in some cases. }
\end{remark}

\noindent{\bf Proof of Theorem \ref{theorem2}.} Define a
$C^2-$function $V: \RR^n_+\rightarrow \RR_+$ by
$$
V(x):=\Pi_{i=1}^nx_i^{p_i}.
$$
Compute
\begin{align}
{\cal L}V(x)&=V(x)p^T(b+Ax)-\frac{1}{2}V(x)x^T\sigma^T[{\rm diag}(p_1,\ldots, p_n)-pp^T]\sigma x\nonumber \\
&\quad+V(x)\int_{\mathbb{Y}}\bigg[ \Pi_{i=1}^n (1+H_i(x, u))^{p_i}-1- \sum_{i=1}^np_iH_i(x, u)\bigg]\lambda(d u)\nonumber \\
&=V(x)p^T(b+Ax)-\frac{1}{2}V(x)x^T\sigma^T[{\rm diag}(p_1,\ldots, p_n)-pp^T]\sigma x\nonumber \\
&\quad+V(x)\sum_{i=1}^n\int_{\mathbb{Y}}\bigg[ (1+H_i(x, u))^{p_i}-1- p_iH_i(x, u)\bigg]\lambda(d u)\nonumber \\
&\quad+V(x)\int_{\mathbb{Y}}\bigg[ \Pi_{i=1}^n (1+H_i(x,
u))^{p_i}-\sum_{i=1}^n(1+H_i(x, u))^{p_i}+n-1\bigg]\lambda(du).
\end{align}
Noting conditions \eqref{eq11} and \eqref{eq11a}, we derive that
there exist positive constants $C_1$ and $C_2$ such that
$$
{\cal L}V(x)\le V(x)(C_1-C_2|x|^\alpha).
$$
For each $k>|\bar{x}|$ define a stopping time
$$
\tau_{k}:=\inf\{t\ge 0: |X(t)|>k\}.
$$
By Theorem \ref{general2} $\tau_k<\infty$ and $\tau_k\rightarrow
\infty$ as $k \rightarrow \infty$ almost surely. Using the It\^o
formula we obtain
\begin{align*}
\mathbb{E}V(X(t\wedge \tau_k))&=V(X_0)+\mathbb{E}\int_0^{t\wedge \tau_k}{\cal L}V(X(s))ds\\
&\le V(X_0)+C_1\mathbb{E}\int_0^{t\wedge \tau_k} V(X(s))ds.
\end{align*}
Hence applying the well-known Gronwall inequality and letting $k\rightarrow \infty$ gives
$$
\mathbb{E}V(X(t))\le V(X_0)e^{C_1t},
$$
and the required assertion follows.

\section{Asymptotic Pathwise Estimation}

In the last section we have discussed how the solutions vary in
$\mathbb{R}^n_+$ in probability or in moment. In this section we
examine pathwise properties of the solutions. To discuss the
pathwise properties of Eq. \eqref{eq1}, we cite the following
exponential martingale inequality with jumps, , e.g., \cite[Theorem
5.2.9, p291]{a09}.

\begin{lemma}\label{exponential martingale}
{\rm Assume that $g:[0,\infty)\rightarrow \mathbb{R}$ and
$h:[0,\infty)\times \mathbb{Y}\rightarrow \mathbb{R}$ are both
predictable $\mathcal {F}_t$-adapted processes such that for any
$T>0$
\begin{equation*}
\int_0^T|g(t)|^2dt< \infty\ \ \mbox{ a.s. and }
\int_0^T\int_{\mathbb{Y}}|h(t,u)|^2\lambda(du)dt<\infty  \ \ \mbox{
a.s.}
\end{equation*}
Then for any constants $\alpha,\beta>0$
\begin{equation*}
\begin{split}
\mathbb{P}\Big\{\sup\limits_{0\leq t\leq
T}\Big[&\int_0^tg(s)dW(s)-\frac{\alpha}{2}\int_0^t|g(s)|^2ds+\int_0^t\int_{\mathbb{Y}}h(s,u)\tilde{N}(ds,du)\\
&-\frac{1}{\alpha}\int_0^t\int_{\mathbb{Y}}[e^{\alpha
h(s,u)}-1-\alpha h(s,u)]\lambda(du)ds\Big]>\beta\Big\}\leq
e^{-\alpha\beta}.
\end{split}
\end{equation*}
}
\end{lemma}

\begin{theorem}\label{exponent}
{\rm Let conditions of Theorem \ref{general2} hold. Assume further
that there exists constant $\theta\in(0,\alpha]$ such that
\begin{equation}\label{eq90}
\int_{\mathbb{Y}}[(\ln Q(x,u))^2+Q(x,u)]\lambda(du)\leq
K|x|^{\theta}+\mbox{o}(|x|^{\theta}),
\end{equation}
where, for $p\in(0,1)$ and $x\in\mathbb{R}_+$
\begin{equation}
Q(x,u):=\sum\limits_{i=1}^n(1+H_i(x,u))^px_i^p\Big/\sum\limits_{i=1}^nx_i^p.
\end{equation}
 There exists $K>0$, independent of initial
value $\bar{x}\in \mathbb{R}^n_+$, such that the solution
$X(t),t\geq0$, of Eq. \eqref{eq1} has the property
\begin{equation}\label{eq16}
\limsup_{t\rightarrow\infty}\frac{\ln(|X(t)|)}{\ln t}\leq K, \ \ \
\mbox{a.s.}
\end{equation}
}
\end{theorem}

\begin{remark}
{\rm By inequality \eqref{eq109}, it is easy to see that condition
\eqref{eq90} holds for $H_i$ in Example \ref{example}. }
\end{remark}

In what follows we complete the argument of Theorem \ref{exponent}.

\noindent{\bf Proof of Theorem \ref{exponent}} Note by Theorem
\ref{general2}  that Eq. \eqref{eq1} has a unique global positive
solution for any initial value $\bar{x}\in\mathbb{R}^n_+$. Let
$V(x), K_1(x,p)$ be defined by \eqref{eq107} and \eqref{eq112},
respectively, and
$Z(x):=\frac{1}{V(x)}\sum_{i=1}^npx_i^p\sum_{j=1}^n\sigma_{ij}x_j$
for $x\in\mathbb{R}^n_+$ and $p\in(0,1)$.  Applying the It\^o
formula to $e^t \ln V(x)$ yields
\begin{equation}\label{eq73}
\begin{split}
e^t\ln V(X(t))&=\ln V(\bar{x})+\int_0^te^s\Big[\ln V(X(s))+\frac{1}{V(X(s))}K_1(X(s),p)-\frac{1}{2}Z^2(X(s))\\
&\quad+\int_{\mathbb{Y}}\Big(\ln
Q(X(s),u)-\frac{p}{V(X(s))}\sum\limits_{i=1}^nX_i^p(s)H_i(X_s,u)\Big)\lambda(du)\Big]ds\\
&\quad+\int_0^te^sZ(X(s))dW(s)+\int_0^t\int_{\mathbb{Y}}e^s\ln
Q(X(s^-),u)\tilde{N}(ds,du).
\end{split}
\end{equation}
By virtue of  Lemma \ref{exponential martingale}, for any
$\alpha,\beta,T>0$ we have
\begin{equation*}
\begin{split}
\mathbb{P}\Big\{\omega:&\sup\limits_{0\leq t\leq
T}\Big[\int_0^te^sZ(X(s))dW(s)-\frac{\alpha}{2}\int_0^te^{2s}Z^2(X(s))ds\\
&\quad+\int_0^t\int_{\mathbb{Y}}e^s\ln
Q(X(s^-),u)\tilde{N}(ds,du)\\
&-\frac{1}{\alpha}\int_0^t\int_{\mathbb{Y}}\Big(Q^{\alpha
e^s}(X(s),u) -1-\alpha e^s\ln
Q(X(s),u)\Big)\lambda(du)ds\Big]\geq\beta\Big\}\leq
e^{-\alpha\beta}.
\end{split}
\end{equation*}
Choose $T=k, \alpha=\epsilon e^{-k}$ and $\beta=\frac{2 e^{k}\ln
k}{\epsilon}$, where $k\in\mathbb{N}, 0<\epsilon<\frac{1}{2}$, in
the above equation. Since $\sum_{k=1}^\infty k^{-2}<\infty$, we can
deduce from the Borel-Cantalli lemma that there exists an
$\Omega_0\subseteq\Omega$ with $\mathbb{P}(\Omega_0)=1$ such that
for any $\omega\in\Omega_0$ we can find an integer $k_0(\omega)>0$
such that
\begin{equation}\label{eq74}
\begin{split}
&\int_0^te^sZ(X(s))dW(s)+\int_0^t\int_{\mathbb{Y}}e^s\ln
Q(X(s^-),u)\tilde{N}(ds,du)\\
&\leq\frac{2 e^{k}\ln k}{\epsilon}+\frac{\epsilon e^{-k}}{2}\int_0^te^{2s}Z^2(X(s))ds\\
&\quad+\frac{1}{\epsilon
e^{-k}}\int_0^t\int_{\mathbb{Y}}\Big(Q^{\epsilon e^{s-k}}(X(s),u)
-1-\epsilon e^{s-k}\ln Q(X(s),u)\Big)\lambda(du)ds
\end{split}
\end{equation}
whenever $0\leq t\leq k$ and $ k\geq k_0(\omega)$. Hence, for any
$\omega\in\Omega_0$, $0\leq t\leq k$ and $ k\geq k_0(\omega)$
\begin{equation*}
\begin{split}
\ln V(X(t))&\leq e^{-t}\ln V(\bar{x})+\frac{2 e^{k-t}\ln k}{\epsilon }\\
&\quad+\int_0^te^{s-t}\Big[\ln V(X(s))+\frac{1}{V(X(s))}K_1(X(s),p)-\frac{1-\epsilon}{2}Z^2(X(s))\Big]ds\\
&\quad+\int_0^te^{s-t}\int_{\mathbb{Y}}\Big(\ln
Q(X(s),u)-\frac{p}{V(X(s))}\sum\limits_{i=1}^nX_i^p(s)H_i(X(s),u)\Big)\lambda(du)ds\\
&\quad+\frac{1}{\epsilon
e^{t-k}}\int_0^t\int_{\mathbb{Y}}\Big(Q^{\epsilon e^{s-k}}(X(s),u)
-1-\epsilon e^{s-k}\ln
Q(X(s),u)\Big)\lambda(du)ds\\
&:=J_1(t)+J_2(t)+J_3(t)+J_4(t).
\end{split}
\end{equation*}
For any $x\in\mathbb{R}^n_+$ and $u\in\mathbb{Y}$, compute  that
\begin{equation}\label{eq117}
\begin{split}
&\ln
Q(x,u)-\frac{p}{V(x)}\sum\limits_{i=1}^nx_i^pH_i(x,u)\\
&=\log
Q(x,u)-Q(x,u)+1+Q(x,u)-\frac{p}{V(x)}\sum\limits_{i=1}^nx_i^pH_i(x,u)-1\\
&\leq
Q(x,u)-\frac{p}{V(x)}\sum\limits_{i=1}^nx_i^pH_i(x,u)-1\\
&=\frac{1}{V(x)}\sum\limits_{i=1}^n[(1+H_i(x,u))^p-1-pH_i(x,u)]x_i^p,
\end{split}
\end{equation}
where in the second  step we used the inequality \eqref{eq110}.  In
the light of a Taylor's series expansion, for
$\epsilon\in(0,\frac{1}{2}]$, $x\in\mathbb{R}^n_+, u\in\mathbb{Y}$
and $s\leq k$
\begin{equation*}
Q^{\epsilon e^{s-k}}(x,u)=1+\epsilon e^{s-k}\ln
Q(x,u)+\frac{\epsilon^2e^{2(s-k)}}{2}(\ln Q(x,u))^2Q^\xi(x,u),
\end{equation*}
where $\xi$ lies between $0$ and $\epsilon$. Thus
\begin{equation*}
J_4(t)=\int_0^t\int_{\mathbb{Y}}\frac{\epsilon e^{2s-t-k}}{2}(\ln
Q(x,u))^2Q^\xi(x,u)\lambda(du)ds.
\end{equation*}
Note that for any $\omega\in\Omega_0$, $t\leq k$ and $ k\geq
k_0(\omega)$
\begin{equation*}
\begin{split}
J_4(t)&=\int_0^t\int_{0<Q(x,u)<1}\frac{\epsilon e^{2s-t-k}}{2}(\ln
Q(x,u))^2Q^{\xi}(x,u)\lambda(du)ds\\
&\quad+\int_0^t\int_{Q(x,u)\geq1}\frac{\epsilon e^{2s-t-k}}{2}(\ln
Q(x,u))^2Q^{\xi}(x,u)\lambda(du)ds\\
&=:\Gamma_1(t)+\Gamma_2(t).
\end{split}
\end{equation*}
For $0<Q(x,u)<1$ and $0\leq\xi\leq\epsilon\leq\frac{1}{2}$, we have
$Q^{\xi}(x,u)\leq1$. Hence
\begin{equation*}
\begin{split}
\Gamma_1(t)&\leq\int_0^t\int_{0<Q(x,u)\leq1}\frac{\epsilon
e^{2s-t-k}}{2}(\ln
Q(x,u))^2\lambda(du)ds\\
&\leq\int_0^t\int_{\mathbb{Y}}\frac{\epsilon e^{2s-t-k}}{2}(\ln
Q(x,u))^2\lambda(du)ds.
\end{split}
\end{equation*}
On the other hand, recalling the fundamental inequality
\begin{equation*}
\ln x\leq 4(x^{\frac{1}{4}}-1) \mbox{ for } x\geq1,
\end{equation*}
and observing $Q^\xi(x,u)\leq Q^{\frac{1}{2}}(x,u)$ for
$Q(x,u)\geq1$ and $0\leq\xi\leq\epsilon\leq\frac{1}{2}$, we have
\begin{equation*}
\begin{split}
\Gamma_2(t)&\leq16\int_0^t\int_{Q(x,u)\geq1}\frac{\epsilon
e^{2s-t-k}}{2}Q(x,u)\lambda(du)ds\\
&\leq16\int_0^t\int_{\mathbb{Y}}\frac{\epsilon
e^{2s-t-k}}{2}Q(x,u)\lambda(du)ds.
\end{split}
\end{equation*}
Consequently, for any $\omega\in\Omega_0$, $t\leq k$ and $ k\geq
k_0(\omega)$
\begin{equation*}
J_4(t)\leq\frac{\epsilon}{2}\int_0^t\int_{\mathbb{Y}} e^{s-t}[(\ln
Q(x,u))^2+16Q(x,u)]\lambda(du)ds.
\end{equation*}
Then, for any $\omega\in\Omega_0$, $0\leq t\leq k$ and $ k\geq
k_0(\omega)$, we have
\begin{equation*}
\begin{split}
J_2(t)+J_3(t)+J_4(t)&\leq\int_0^te^{s-t}\Big[\ln
V(X(s))+\frac{1}{V(X(s))}K_1(X(s),p)-\frac{1-\epsilon}{2}Z^2(X(s))\\
&\quad+\int_{\mathbb{Y}}\sum\limits_{i=1}^n\frac{X_i(s)^p}{V(X(s))}[(1+H_i(X(s),u))^p-1-pH_i(X(s),u)]\lambda(du)\\
&\quad+\frac{\epsilon}{2}\int_{\mathbb{Y}}[(\ln
Q(X(s),u))^2+16Q(X(s),u)]\lambda(du)\Big]ds\\
&:=\int_0^te^{s-t}[M_1(X(s),p)+M_2(X(s),p)+M_3(X(s),p)]ds.
\end{split}
\end{equation*}
By inequality \eqref{eq109} and inequality \eqref{eq110}, for
$x\in\mathbb{R}^n_+$ and $p\in(0,1)$ we obtain
\begin{equation*}
M_1(x,p)\leq K|x|^2+\mbox{o}(|x|^2).
\end{equation*}
Moreover, again by inequality \eqref{eq109}, together with
\eqref{eq11} and \eqref{eq90}, we can deduce that for
$x\in\mathbb{R}^n_+$ and $\alpha>2$
\begin{equation*}
M_2(x,p)+M_3(x,p)\leq -K|x|^\alpha+\mbox{o}(|x|^\alpha).
\end{equation*}
Noting that $M_1(x,p)+M_2(x,p)+M_3(x,p)$ is bounded by a polynomial
with the negative leading coefficient, we arrive at
\begin{equation*}
J_2(t)+J_3(t)+J_4(t)\leq K\int_0^te^{s-t}ds=K(1-e^{-t}).
\end{equation*}
Thus, for any $\omega\in\Omega_0, 0\leq t\leq k$ and $ k\geq
k_0(\omega)$
\begin{equation*}
\begin{split}
\ln V(X(t))&\leq e^{-t}\ln V(\bar{x})+\frac{2 e^{k-t}\ln k}{\epsilon
}+K.
\end{split}
\end{equation*}
In particular, for $\omega\in\Omega_0, k-1\leq t\leq k$ and $k\geq
k_0(\omega)+1$, we have
\begin{equation*}
\begin{split}
\frac{\ln(|X(t)|^p)}{\ln t}\leq\frac{\ln V(X(t))}{\ln t}&\leq
\frac{1}{\ln(k-1)}\left[e^{-t}\ln V(\bar{x})+K\right]+\frac{2 e\ln
k}{\epsilon\ln(k-1) }.
\end{split}
\end{equation*}
This implies that
\begin{equation*}
\limsup_{t\rightarrow\infty}\frac{\ln(|X(t)|)}{\ln t}\leq\frac{2
e}{p\epsilon }.
\end{equation*}
The desired assertion then follows by letting
$\epsilon\uparrow\frac{1}{2}$.

Noting the limit $\lim_{t\rightarrow\infty}\frac{\ln t}{t}=0$, we
can easily deduce from Theorem \ref{exponent} that the sample
Lyapunov exponent of Eq. \eqref{eq1} is less or equal to zero, which
is stated as the following corollary.
\begin{corollary}
{\rm Under conditions of Theorem \ref{exponent},
\begin{equation*}
\limsup_{t\rightarrow\infty}\frac{\ln(|X(t)|)}{t}\leq 0, \ \ \
\mbox{a.s.}
\end{equation*}
}
\end{corollary}

\begin{remark}
{\rm Our theories developed can also be applied to discuss
stochastic functional Kolmogorov-type population dynamics with jumps
and switching-diffusion ecosystems with jumps, respectively, which
will be reported in separated papers. }
\end{remark}

\end{document}